\documentstyle[12pt]{article}
\begin{document}
\large
\begin{center}
{\bf EXPONENTIAL BOUNDS IN THE LAW OF ITERATED }
\vspace{3mm}
  {\bf LOGARITHM FOR MARTINGALES }\\
\vspace{3mm}

{\bf E. Ostrovsky, L.Sirota.}\\

\vspace{3mm}
Department of Mathematics, Ben - Gurion University, 84105, Beer - Sheva, Israel.\\
e-mail: \ galo@list.ru \\
Department of Mathematics, Bar - Ilan University, Ramat - Gan, 59200,Israel.\\
e - mail: \ sirota@zahav.net.il\\

\end{center}
\vspace{3mm}
{\bf Abstract.}  In this paper non-asymptotic exponential
estimates are derived for  tail of maximum martingale distribution by naturally
norming in the spirit of the classical Law of Iterated Logarithm.  


\vspace{3mm}
\begin{center}
{\bf 1. Introduction. Notations. Statement of problem.} \par
\end{center}
\vspace{3mm}

 Let $ (\Omega,F,{\bf P} ) $ be a probability space, $  \Omega = \{\omega\}, \
(S(n), F(n)), \ n = 1,2, \ldots  $ being a centered: $ {\bf E} S(n) = 0 \ $ non-trivial:

$$
\forall n \ \Rightarrow \sigma(n) =  [ {\bf Var} \ (S(n))]^{1/2} \in (0, \infty)
$$
 martingale: $ {\bf E} S(n+1)/F(n) = S(n) $ relatively some filtration $ F(n).$
Let also $ v(n) $ be a {\it deterministic} positive monotonically increasing sequence, $ A(k) $ be a
deterministic positive strong monotonically increasing {\it integer}
sequence $ A(k), \ k = 1,2, \ldots $ such that
 $ A(1) = 1, \ B(k) \stackrel{def}{=} A(k+1) - 1  \ge A(k) + 1. \  $ Introduce the
{\it} partition of integer semi-axis $ Z_+ = [1,2,\ldots ) \ R = \{A(k), B(k) \}: $

$$
Z_+ =  \cup_{k=1}^{\infty}[A(k), B(k)]  = \cup_{k=1}^{\infty} [A(k), A(k+1) - 1]
$$
and denote the {\it set} of all these partitions by $ T: \ T = \{ R \}.$ \par
 Let us introduce the following probability $ W(u): $

$$
W(u) = W(v; u) \stackrel{def}{=} {\bf P} \left( \sup_n \frac{S(n)}{\sigma(n) \ v(n)}  > u \right), \eqno(1)
$$

and analogously set

$$
W_+(u) = W_+(v; u) \stackrel{def}{=} {\bf P} \left( \sup_n \frac{|S(n)|}{\sigma(n) \ v(n)}  > u \right).
$$

{\bf Our goal is obtaining the exponential decreasing estimation for $ W(u),
 \ W_+(v,u) $ for sufficiently greatest values $ u, $ for example,
 $ u \ge 2. $ } \par
 In the case when $ S(n) = \sum_{i=1}^n \xi(i), $ where $ \{ \xi(i) \} $ are
independent centered r.v. and  $ \sigma \ - $ flow $ \{ F(n) \} $ is the
natural filtration:
$$
F(n) = \sigma \{ \xi(i), \ i = 1,2, \ldots,n  \}
$$
with the classical norming  $ v(n) = (\log (\log(n + 3))^{1/2} $ the estimation
for  $  P(u) $ was obtained in [1], see also [2], p.62 - 66. Our result
may be considered as some addition to the classical Law of Iterated Logarithm
(LIL) for martingales, i.e. of the view

$$
\overline{\lim}_{n \to \infty} \ |S(n)| / (\sigma(n) \ v(n) ) = \eta(\omega) < \infty \ {\bf a.e.}, \eqno(2)
$$
see [3], p.115-127  and references there. \par
 It is clear that if the conclusion (2) is satisfied, then the bound for
$ P(u) $ is not trivial, i.e. $ u \to \infty \ \Rightarrow P(u) \to 0. $ \par

\vspace{2mm}

{\bf 2. Result.} \par
 In order to formulate our result, we need to introduce some another notations
and conditions. Let $ \phi = \phi(\lambda), \lambda \in (-\lambda_0, \lambda_0),  \ \lambda_0 = const \in (0, \infty] $ be some even taking positive values
for positive arguments strong convex twice continuous differentiable function,
such that
$$
 \phi(0) = 0, \ \lim_{\lambda \to \lambda_0} \phi(\lambda)/\lambda = \infty.
\eqno(3)
$$
 The set of all these function we denote $ \Phi; \ \Phi = \{ \phi(\cdot) \}. $
We say that the {\it centered} random variable (r.v) $ \xi = \xi(\omega) $
belongs to the space $ B(\varphi), $ if there exists some non-negative constant
$ \tau  \ge 0 $  such that

$$
\forall \lambda \in (-\lambda_0, \lambda_0)  \ \Rightarrow
{\bf E} \exp(\lambda \xi) \le \phi(\lambda \ \tau). \eqno(4).
$$
 The minimal value $ \tau $ satisfying (4) is called the $ B(\phi) \ $ norm
of the variable $ \xi, $ write
 $$
 ||\xi||B(\phi) = \inf \{ \tau, \ \tau > 0: \ \forall \lambda \ \Rightarrow  {\bf E}\exp(\lambda \xi) \le \exp(\phi(\lambda \ \tau))  \}.
 $$
 This spaces are very convenient for the investigation of the r.v. having a
exponential decreasing
 tail of distribution, for instance, for investigation of
limit theorem,  exponential bounds of distribution for sums of random variables,
non-asymptotical properties, problem of continuous  of random fields etc.\par

 The space $ B(\phi) $ relative to the norm $ || \cdot ||B(\phi) $  is a
Banach space which is isomorphic to subspace consisted on all the centered
variables of Orlich's space $ (\Omega,F,{\bf P}), N(\cdot) $ with $ N \ - $ function
$$
N(u) = \exp(\phi^*(u)) - 1, \ \phi^*(u) = \sup_{\lambda} (\lambda u -
\phi(\lambda)).
$$
 The transform $ \phi \to \phi^* $ is called Young-Fenchel transform. The proof of considered assertion used the properties of saddle-point method and theorem
 of Fenchel-Moraux:
$$
\phi^{**} = \phi.
$$

 The next facts about the $ B(\phi) $ spaces are proved in [2, p. 19 - 40], [4]:

$$
{\bf 1.} \ \xi \in B(\phi) \Leftrightarrow {\bf E } \xi = 0, \ {\bf and} \ \exists C = const > 0,
$$

$$
U(\xi,x) \le \exp(-\phi^*(Cx)), x \ge 0,
$$
where $ U(\xi,x)$ denotes as usually the tail of distribution of a r.v. $ \xi: $

$$
U(\xi,x) = \max \left( {\bf P}(\xi > x), \ {\bf P}(\xi < - x) \right),
\ x \ge 0, \eqno(5)
$$
and this estimation (5) is in general case asymptotically exact. \par
 Here and further $ C, C_j, C(i) $ will denote the non-essentially positive
finite "constructive" constants.\par
 More exactly, if $ \lambda_0 = \infty, $ then the following implication holds:

$$
\lim_{\lambda \to \infty} \phi^{-1}(\log {\bf E} \exp(\lambda \xi))/\lambda =
K \in (0, \infty)
$$
if and only if

$$
\lim_{x \to \infty}  (\phi^*)^{-1}( |\log U(\xi,x)| )/x = 1/K.
$$
 Here and further $ f^{-1}(\cdot) $ denotes the inverse function to the
function $ f $ on the left-side half-line  $ (C, \infty). $ \par

{\bf 2.} Define $ \psi(p) = p/\phi^{-1}(p), \ p \ge 2.  $
 Let us introduce the new norm on the set of r.v. defined in our probability space by the following way: the space $ G(\psi) $ consist, by definition, on all the centered r.v. with finite norm

$$
||\xi||G(\psi) \stackrel{def}{=} \sup_{p \ge 2} |\xi|_p/\psi(p), \ |\xi|_p =
{\bf E}^{1/p} |\xi|^p. \eqno(6)
$$
 It is proved that the spaces $ B(\phi) $ and $ G(\psi) $ coincides:$ B(\phi) =
G(\psi) $ (set equality)  and both
the norm $ ||\cdot||B(\phi) $ and $ ||\cdot|| $ are equivalent: $ \exists C_1 =
C_1(\phi), C_2 =  C_2(\phi) = const \in (0,\infty), \ \forall \xi \in B(\phi) $

$$
||\xi||G(\psi) \le C_1 \ ||\xi||B(\phi) \le C_2 \ ||\xi||G(\psi).
$$

{\bf 3.} The definition (6) is correct  still for the non-centered random
variables $ \xi.$ If for some non-zero r.v. $ \xi \ $ we have $ ||\xi||G(\psi) < \infty, $ then for all positive values $ u $

$$
{\bf P}(|\xi| > u) \le 2 \ \exp \left( - u/(C_3 \ ||\xi||G(\psi))  \right).
\eqno(7)
$$
and conversely if a r.v. $ \xi $ satisfies (7), then $ ||\xi||G(\psi) <
\infty. $ \par
{\sc We suppose in this article that there exists the function $ \phi \in \Phi $
such that}
$$
\sup_n [ ||S(n)||B(\phi)/\sigma(n)] < \infty,
$$
or equally for all non-negative values $  x $

$$
\sup_n \max \left[{\bf P}\left( \frac{S(n)}{\sigma(n)} > x \right), \ {\bf P} \left( \frac{S(n)}{\sigma(n)} < - x \right) \right] \le \exp \left(-\phi^*(x/C) \right).  \eqno(8)
$$
 The function $ \phi(\cdot) $ may be constructive" introduced by the formula
$$
\phi(\lambda) = \log \sup_n {\bf E} \exp(\lambda S(n)/\sigma(n)),
$$
 if obviously the family of r.v. $ \{S(n)/\sigma(n)\} $ satisfies the
{\it uniform } Kramer's condition: $ \exists \mu \in (0,\infty), \ 
\forall x > 0 \ \Rightarrow $
$$
 \sup_n U(S(n)/\sigma(n), \ x) \le \exp(-\mu \ x).
$$

 There are many examples of martingales satisfying  the condition (8) in the article [5]; in particular, there are many examples with
$$
\phi^*(x) = x^r \ L(x),  \ r = const > 0, \eqno(9)
$$

$$
n^{\gamma} \ M_1(n) \le \sigma(n) \le n^{\gamma} \ M_2(n), \ \gamma =
const > 0, \eqno(10)
$$
where $ L(x), M_1(n), M_2(n) $ are some positive continuous {\it slowly
varying} as $ x \to \infty $ or correspondently as $ n \to \infty $
functions. \par

 Let us denote for some partition $ R = \{ A(k), B(k) \}  $

$$
 Q(k; R,v,u) = \exp \left(-\phi^*(u \sigma(A(k)) \ v(A(k))/\sigma(B(k) )  \right),
$$

$$
Q(R,v,u) = \sum_{k=1}^{\infty} Q(k; R,v,u). \eqno(11)
$$

{\bf Theorem.} {\it Under  our conditions and for some finite } $ C = C(\phi) $

$$
W(v;u) \le \inf_{R \in T} Q(R,v,Cu), \eqno(12)
$$
{\it and analogous estimation is true for the probability } $ W_+(v,u). $ \par
 {\bf Proof.} Let $ Z_+ = \cup_k [A(k), B(k)], \ B(k) = A(k+1) - 1 $ be arbitrary
partition, $ R = \{A(k), B(k) \}  \in T.$  Denote  $ E(k) = [A(k), B(k)].$ We see:

$$
W(v;u) \le \sum_{k=1}^{\infty}W(k; v,u), \ W(k; v,u) \stackrel{def}{=}
$$
$$
{\bf P } \left( \max_{n \in E(k)} (S(n)/(\sigma(n) \ v(n)) > u \right). \eqno(13)
$$
 Let us estimate the probability $ W(k; v,u). $ We obtain:

$$
W(k; v,u) \le {\bf P} \left(\max_{n \in E(k)} S(n) > u \ \sigma(A(k)) \ v(A(k))
/\sigma(B(k)) \right),
$$
as long as both the functions $ \sigma(\cdot) $ and $ v(\cdot) $
are  monotonically increasing.  \par
 Further we use the Doob's inequality and properties of $ B(\phi) $ spaces.
It follows from Doob's inequality
$$
| \max_{n \in E(k)} S_n|_p \le C \ \sigma(B(k)) \cdot (p/\phi^{-1}(p)) \cdot
(p/(p-1)) \le
$$
$$
2 \ C \ \sigma(B(k)) \cdot (p/\phi^{-1}(p))
$$
as long as $ p \ge 2.$  Therefore $ W(k; v,u) \le $

$$
  \exp \left(-\phi^*(C u \ \sigma(A(k)) \ v(A(k))/\sigma(B(k) )  \right) = Q(k; R,v,Cu).
$$
We obtain alter summation

$$
W(v; u) \le Q(R,v,Cu).
$$
 Since the partition $ R $ is arbitrary, we get to the  demanded inequality
(12). \par
 The probability $ W_+(v;u) $ is estimated analogously, as long as
$ (- S(n), F(n)) $ is again the martingale with at the same function
$ \phi(\cdot). $ \par
 Note that we can ground our theorem from the Kolmogorov's inequality for martingales.\par

{\bf 3.  Examples.} Let us consider some examples in order to show the exactness
of our theorem. \par
{\bf A.} Let $ \eta $ be a symmetrically distributed  r.v. with the tail of
distribution of a view:

$$
   {\bf P}(\eta > x)  = \exp \left( - \phi^*(x) \right),
$$
$ x \ge 0, \ \phi \in \Phi; $ and let $ \{ \xi(i) \} $ be an independent copies
of $ \eta.$ Then $ ||\eta||B(\phi) = C_5 \in (0,\infty), \ \beta^2 \stackrel{def}{= }
{\bf Var} (\eta) \in (0,\infty). $ \par
 Let us consider the martingale $ (S(n), F(n)), $ where
$$
S(n) = \sum_{k=1}^n  2^{-k} \ \xi(k)
$$
relative the natural filtration $ \{ F(n) \}. $ It follows from the triangle
inequality for the $ B(\phi) $ norm that
$$
\sup_n ||S(n)||B(\phi) \le \sum_{k=1}^{\infty} 2^{-k} \ ||\xi(k)||B(\phi) =
C_5 < \infty,
$$

$$
0.25 \ \beta^2 \le \sigma^2(n) \le \beta^2;
$$

therefore
$$
\exp \left(- \phi^*(C_6 x) \right) \le \sup_n {\bf P}( S(n)/\sigma(n) > x) \le
$$
$$
\exp \left(- \phi^*(C_7 x) \right),
$$

 $ 0 < C_7 < C_6 < \infty $ (the low bound is trivial).\par
 Moreover, it is possible to prove that

$$
\inf_n {\bf P}( S(n) > x) \ge \exp \left(- \phi^*(C_8 x) \right).
$$

{\bf B.} Assume here that the martingale $ ( S(n), F(n) ) $ satisfies the
conditions (9) and (10). Let us choose

$$
v(n) = v_r(n) = [ \log(\log(n+3))]^{1/r},
$$
or equally

$$
v(n) = v_r(n) = [ \log(\log( \sigma(n)+3))]^{1/r},
$$

then  we obtain after some calculation on the basis of our theorem, choosing
the partition $ R = \{ [A(k), A(k+1) - 1] \} $ such that:
$$
A(k) = Q^{k-1},
$$
where $  Q = 3 $ or $ Q = 4 $ etc.:
$$
{\bf P}\left(\sup_n \frac{S(n)}{\sigma(n) \ v_r(n) } > x \right) \le
\exp \left[- C \  x^r \ L(x) \right], x > 0. \eqno(14)
$$
 Moreover, if the martingale $ (S(n), F(n) $ satisfies the conditions (8),
 (9) and (10), then with probability one
$$
\overline{\lim}_{n \to \infty} \frac{S(n)}{\sigma(n) \ v_r(n)} \le C,
$$
where the constant $ C $  is defined in (8);  and
the last inequality is exact, e.g., for the martingales considered
in the next  section {\bf C.} \par

{\bf C.} Let us show the exactness of the estimation (14). Consider the
so-called Rademacher sequence $ \{ \epsilon(i) \}, \ i = 1,2,\ldots; $ i.e. where
$ \{ \epsilon(i) \} $ are independent and $ {\bf P}(\epsilon(i) = 1) =
{\bf P}(\epsilon(i) = - 1) = 0.5. $ \par
  It is known that that the r. v. $ \{ \epsilon(i) \} $ belongs to the
$ B(\phi_2) $ space with corresponding function
$$
\phi_2(\lambda) = 0.5 \ \lambda^2, \ \lambda \in  (-\infty, \infty).
$$

Denote for $ d = 1,2,3,\ldots \ S(n) = S_d(n) = $

$$
 \sum \sum \ldots \sum_{1 \le i(1)< i(2) \ldots < i(d) \le n}
 \epsilon(i(1)) \ \epsilon(i(2)) \ \epsilon(i(3)) \ldots \ \epsilon(i(d))
$$
under natural filtration $ F(n). $  It is easy to verify that $ (S(n), F(n)) $
is a martingale and that
 $$
        0 < C_1 \le \sigma^2(n)/n^d \le C_2 < \infty.
 $$

It follows from our theorem that

$$
{\bf P} \left( \sup_{n} \frac{S(n)}{ (n \ \log(\log(n + 3)))^{d/2} } > u \right) < \exp \left(- C u^{2/d} \right),
$$
and as it is proved in [5]

$$
\exp \left[- C_3 \  x^{2/d} \right] \le
$$

\vspace{3mm}

$$
\sup_n {\bf P}\left(\frac{|S(n|}{\sigma(n)} > x \right) \le
\exp \left[- C_4 \  x^{2/d} \right], x > 0,
$$

i.e. in the considered case $ r = 2/d. $ \par
 We prove in addition that
$$
{\bf P} \left( \overline {\lim}_{n \to \infty}
 \frac{S(n)}{(n \ \log ( \log(n + 3) ))^{d/2} } > 0 \right) > 0. \eqno(15)
$$
 Ii is enough to consider only the  case $ d = 2, $  i.e. when
 $$
 S(n) = \sum \sum_{1 \le i < j \le n} \epsilon(i) \ \epsilon(j).
 $$

 We observe that
 $$
 2 \ S(n) = \left(\sum_{k=1}^n \epsilon(k) \right)^2 - \sum_{m=1}^n  (\epsilon(m))^2 \stackrel{def}{=} \Sigma_1(n) - \Sigma_2(n).
 $$

 From the classical LIL on the form belonging to Hartman-Wintner it follows
 that there exist a finite
 non-trivial non-negative  random  variables $ \theta_1, \ \theta_2 $ for which

$$
|\Sigma_2(n)| \le n + \theta_2 \sqrt{n \ \log(\log(n + 3 )) } \eqno(16)
$$

and

$$
\Sigma_1(n_m) \ge \theta_1 \ n_m \ \log(\log(n_m + 3)) \eqno(17)
$$

for some (random) integer positive subsequence  $ n_m, \ n_m \to \infty $ as
$ m \to \infty.$ \par
 The proposition (15) it follows immediately from (16) and (17). \par

 More exactly, by means of considered method may be proved the following
relation:
$$
\overline {\lim}_{n \to \infty} \frac{S(n)}{(n \ \log ( \log(n + 3) ))^{d/2} }
 \stackrel{a.e}{=} \frac{2^{d/2}}{d!}.
$$

{\bf 4.} It is easy to prove the non-improvement of the estimation (14).
Namely, let us consider the martingale $ (S(n), F(n))  $ satisfying the
conditions (9) and (10) and such that for some $ n_0 = 1,2,3,\ldots $
$$
{\bf P} \left( \frac{S(n_0)}{\sigma(n_0)} > u \right) \ge \exp
\left( - C_9 \ u^r \ L(u) \right);
$$
then
$$
W(v_r; u) \ge {\bf P} \left( \frac{S(n_0)}{\sigma(n_0) \ v_r(n_0)} > u \right)
=
$$
\vspace{3mm}
$$
{\bf P}\left( \frac{S(n_0)}{\sigma(n_0)} > u \ v_r(n_0) \right) \ge
\exp \left( - C_{10} \ u^r \ L(u) \right),
$$
since the function $ L(\cdot) $ is slowly varying. \par

\vspace{3 mm}

{\bf 4.  Concluding remarks.} \\

{\bf 1.} It is evident that only the case when
$$
\lim_{n \to \infty} \sigma(n) = \infty
$$
is interest.\par

{\bf 2.} Instead the norm $ \sigma(n) = |S(n)|_2 $ we can consider some
another rearrangement invariant norm in our probability space, say,  the
$ L_s \ $ norm

$$
\sigma_s(n) = |S(n)|_s, s = const \ge 1
$$
 or some norm in Orlicz's space, $ B(\nu), \ \nu \in \Phi $ norm etc. \par
 But the norm $ \sigma(n) $ is classical and more convenient. For instance, if
$ S(0) \stackrel{def}{=} 0, $ then

$$
\sigma^2(n) = \sum_{k=0}^{n - 1} {\bf Var} (S(k+1) - S(k)).
$$

{\bf 3.} The exponential bounds for tail of distribution in the LIL for
martingales used, for instance, in the non-parametric statistic by adaptive estimations (see [6]).

\newpage

{\bf References }\\

   1. Ostrovsky E.I. 1994. {\it Exponential Bounds  in the Law of Iterated                                   Logarithm
      in Banach Space.} Math. Notes, {\bf 56}, 5, p. 98 - 107. \\
   2. Ostrovsky E.I., 1999. Exponential estimations for Random Fields and
     its applications (in Russian). 1999, Obninsk, Russia, OINPE.\\
   3. Hall P., Heyde C.C.,1980. Martingale Limit Theory and Applications.
      Academic Press, New York.\\
   4. Kozachenko Yu. V., Ostrovsky E.I., 1985. The Banach Spaces of
      random Variables of subgaussian type. Theory of Probab. and Math.
      Stat., (in Russian). Kiev, KSU, v.32, 43 - 57.\\
   5. Ostrovsky E. {\it Bide-side exponential and moment inequalities  for tails
      of distribution of Polynomial Martingales.} Electronic publication, arXiv:
      math.PR/0406532 v.1  Jun. 2004. \\
   6. Ostrovsky E., Zelikov Y. {\it Adaptive Optimal Regression and
      Density Estimations based on Fourier-Legendre Expansion.} Electronic
      Publication, arXiv 0706.0881 [math.ST], 6 Jun. 2007.\\

\newpage

\begin{center}

{\bf Ostrovsky E.}\\
\vspace{4mm}

 Address: Ostrovsky E., ISRAEL, 76521, Rehovot, \ Shkolnik street.
5/8. Tel. (972)-8- 945-16-13.\\
\vspace{4mm}
e - mail: {\bf Galo@list.ru}\\

{\bf Sirota L.}\\
\vspace{4mm}
 Address: , Sirota L., ISRAEL, 84105, Beer-Sheeba,

\vspace{4mm}
e - mail: {\bf sirota@zahav.net.il}\\

\end{center}

\end{document}